\documentclass[11pt]{article}

\input{amssym.def}
\input{amssym}
\usepackage{eufrak}
\usepackage{color}

\newtheorem{theorem}{Theorem}

\newtheorem{proposition}[theorem]{Proposition}
\newtheorem{remark}[theorem]{Remark}

\def\pa{\partial}
\def\h{h}
\def\rh{r_\hbar}
\def\De{\Delta}

\def\Cas{{\rm Cas}}

\def\AA{\cal A}
\def\BB{\cal B}

\def\H{{\bf H}}

\def\hh{\hbar}

\def\span{{\rm span}}

\def\AAA{{\bf A}}

\def\de{\delta}

\def\ot{\otimes}

\def\C{{\Bbb C}}

\def\AA{{\cal A}}

\def\Sym{{\rm Sym\, }}

\def\HH{\mathbf{H}}

\def\lhq{\ifmmode {\cal L}(q,\hbar)\else ${\cal L}(q,\hbar)$\fi}

\def\lqh{\ifmmode {\cal L}(q,\hbar)\else ${\cal L}(q,\hbar)$\fi}

\def\hpa{\hat{\pa}}
\def\hD{\hat{D}}
\def\DD{\mathrm{D}}
\def\thD{\hat{D}^t}

\def\hatt{\hat{\Theta}}

\def\la{{\lambda}}
\def\al{{\alpha}}

\def\be{\begin{equation}}
\def\ee{\end{equation}}

\begin{document}

\makeatletter
\renewcommand{\theequation}{{\thesection}.{\arabic{equation}}}
\@addtoreset{equation}{section} \makeatother

\title{Noncommutative geometry on central extension of $U(u(2))$}

\author{\rule{0pt}{7mm} Dimitri
Gurevich\thanks{gurevich@ihes.fr}\\
{\small\it Universit\'e Polytechnique Hauts-de-France,  LMI}\\
{\small\it F-59313 Valenciennes, France}\\
{\small \it and}\\
{\small \it Interdisciplinary Scientific Center J.-V.Poncelet}\\
{\small\it Moscow 119002, Russian Federation}\\
\rule{0pt}{7mm} Pavel Saponov\thanks{Pavel.Saponov@ihep.ru}\\
{\small\it
National Research University Higher School of Economics,}\\
{\small\it 20 Myasnitskaya Ulitsa, Moscow 101000, Russian Federation}\\
{\small \it and}\\
{\small \it
Institute for High Energy Physics, NRC "Kurchatov Institute"}\\
{\small \it Protvino 142281, Russian Federation}}

\maketitle

\begin{abstract}
In our previous publications we have introduced analogs of partial derivatives on the algebras $U(gl(N))$.  In the current paper  we compare two methods of introducing these
analogs: via  the so-called quantum doubles and by means of  a coalgebraic structure. In the case $N=2$ we extend
the quantum partial derivatives from $U(u(2))$ (the compact form of the algebra $U(gl(2))$) on a bigger algebra, constructed in two steps.
First, we define  the derivatives on a  central extension of this algebra,  then we prolongate them on some elements  of the corresponding skew-field by using the 
Cayley-Hamilton identities for certain matrices with noncommutative entries. Eventual applications of this differential calculus are discussed.
\end{abstract}

\section{Introduction}

Attempts to construct a Noncommutative (NC) geometry, related to Quantum Groups (QG) have been undertaken since  the creation of the QG theory. Initiated in
\cite{W}, this study was developed in \cite{IP1, FP} and other papers. The authors of these papers aimed at
constructing  a differential calculus in which the role of the function algebra on the group $GL(N)$ was played by an RTT algebra and that of the one-sided
vector fields was played by elements of the corresponding Reflection Equation (RE) algebra.

In \cite{GPS2} we have constructed a different version of the calculus in which the role of function algebra was attributed to another copy of the RE algebra. By combining the  generating matrices\footnote{Note that all algebras, we are dealing with, are Quantum Matrix (QM) algebras, i.e. they are introduced via systems of relations imposed onto the entries of some matrices. We call them  {\em the generating matrices} .} of these two copies of the RE algebra, we constructed a new matrix, whose entries $\pa_i^j$  played the role of partial derivatives in entries of
the matrix $L=\|l_i^j\|_{1\leq i,j,\leq N}$, which is the generating matrix for the second copy of the RE algebra.

This scheme can be expressed in terms of a Quantum Double (QD) $(A,B)$, where the  algebra $B$ is generated by  the entries $l_i^j$ of the matrix $L$ and the algebra $A$ is generated by the  entries of the matrix $D=\|\pa_i^j\|_{1\leq i,j,\leq N}$. The action of the $\partial_i^j$ on $l_k^s$ is defined via the so-called permutation map
$$
\sigma:A\ot B\to B\ot A,
$$
and  a counit $\varepsilon:A\to \C$. The ground field is assumed to be $\C$.

The RE algebra, we are dealing with,  is defined by entries of a generating matrix $L$ subject to the following relation
\be 
R\, L_1\, R\, L_1-L_1\,R\, L_1\, R=0, 
\label{RE} 
\ee
where $R$ is a  Hecke symmetry. (The reader is referred to \cite{GPS1} for the terminology.)
Note that  by a linear shift of generators it is possible to transform the RE algebra into a modified RE algebra, which is defined by the system
\be 
R\, L_1\, R\, L_1-L_1\,R\, L_1\, R=R\, L_1-L_1\,R. 
\label{mRE} 
\ee
The latter algebra tends to $U(gl(N))$, provided the Hecke symmetry $R$ tends to the usual flip $P$ as $q\to 1$ (see (\ref{odin})). The entries of the matrix $D$  turn at this limit into operators which are treated to be patrial derivatives in $l_i^j$. Below, we call them {\em the quantum partial derivatives} (QPD). Similarly to the classical case, the QPD commute with each other and thus generate a commutative algebra, but the Leibniz rule for them differs from the classical one. Nevertheless, the new Leibniz rule  can be also expressed  via a coproduct $\De$, acting on this algebra (by contrast with the braided versions of $U(gl(N))$, considered in \cite{GPS2}).

The problem of defining analogs of partial derivatives on the algebra $U(g)$, where $g$ is an arbitrary  Lie algebra,  was considered in  \cite{MS}.
Partial derivatives constructed in that paper  were introduced by the  methods of deformation quantization.
The Leibniz rule for the corresponding partial derivatives differs  from ours  even if $g=gl(N)$.
The  coproduct $\De$, entering the Leibniz rule from  \cite{MS}, being applied to  the partial derivatives, gives rise to infinite series in these derivatives. 
Whereas, our Leibniz rule is very simple and can be presented in a group-like form. We do not know other examples of algebras $U(g)$ admitting  partial derivatives with a similar property.

Our global aim is to extend the action of QPD onto a central extension $\AA$ of the algebra $U(gl(N))$ and furthermore on  the skew-field
$$
{\BB}={\AA}[{\AA}^{-1}]=\{a/b\,|\,a, b\in \AA,\, b\not=0\}.
$$

Once this extension is constructed, we become able to transfer a large number of usual differential operators, defined on the algebra $Sym(gl(N))$ onto the algebra $\BB$.
We call this procedure {\em quantization with non-commutative configuration space}. However, since the  Leibniz rule for our QPG does not admit any analog of the formula for
derivatives of  composed functions,  the results of application of the QPD to elements of $\BB$ become more difficult to compute.

 In the current paper we consider the central extension of the compact form  $U(u(2))$ of the algebra $U(gl(2))$, by adding  the "eigenvalues" of its
generating matrix $L$ to the algebra  $U(u(2))$. By these eigenvalues we mean   the roots of the {\it characteristic polynomial} arising from the Cayley-Hamilton (CH) identity 
for the matrix $L$. Note that
a central extension of the RE algebra  (playing the role of the algebra $A$ in the related {\em Heisenberg  double} $(A,B)$) was used in \cite{IP2} in order to study the dynamics of a 
$q$-isotropic top, introduced in \cite{AF}. By contrast with \cite{IP2}, we construct a central extension of the algebra $B$ in the corresponding quantum double $(A,B)$.
It should be emphasized that the CH identities for the matrices with non-commutative entries is the main tool enabling us to extend the action of the QPD onto some elements of $\BB$.
In the last section we compare this method with the Gelfand-Retakh method based on the use of quasideterminants.

The paper is organized as follows.   In the next section we comparer two methods of defining QPD on the algebra $U(gl(N))$: one is based on using a
quantum double with properly defined permutation relations, the other one is based on constructing the corresponding coproduct defined  on the algebra, generated by the QPD.
In section 3 we consider the case $N=2$ and for the compact form $U(u(2))$ of the algebra $U(gl(2))$  we express the Leibniz rule in terms of a matrix $\hatt$,
composed of the QPD, on which the coproduct acts in a multiplicative or group-like way. In section 4 we define a central extension of the algebra $U(u(2))$ by using the 
eigenvalues of the generating matrix $L$ and prolongate on them the actions of the QPD.  Finally, in section 5 we extend our QPD onto some elements from the corresponding skew-field, 
in particular, those  which could be useful for finding a non-commutative counterpart of the Dirac vector-potential for the magnetic monopole field. 

\medskip

\noindent
{\bf Acknowledgement}\par
\noindent
The work of P.S. was partially supported by the RFBR grant 19-01-00726.

\section{Different methods of defining QPD}

Consider the algebra $U(gl(N)_h)$, where $h$ is a numerical factor, introduced in the Lie bracket\footnote{This procedure enables us to treat the algebra $U(gl(N))$ as the result of a quantization (deformation) of the commutative algebra $\Sym(gl(N))$.}. It is generated by the unit $1_{U(gl(N)_h)}$ and  elements $l_i^j$, $1\leq i,j\leq N$, subject to the following system of relations:
\be
l_i^j\, l_k^s-l_k^s\, l_i^j=h(l_i^s\,\de_k^j-l_k^j\, \de_i^s)\quad \Leftrightarrow \quad L_1\, L_2-L_2\,L_1=h(L_1\,  P_{12}-P_{12}\, L_1), 
\label{odin}
\ee
where  $L=\|l_i^j\|_{1\leq i,j\leq N}$,  $L_1=L\ot I, \,\,L_2=I\ot L$. Hereafter,  $I$ stands for the identity operator or its matrix and $P$ stands for the usual flip or its matrix.

The generating matrix $L$ satisfies the Cayley-Hamilton identity. Namely, there exists a monic polynomial
$$
p(t)=t^N+a_1 t^{N-1}+...+a_{N-1}\, t+a_N,\qquad a_k\in Z(U(gl(N)_h)),
$$
called {\em characteristic}, such that $p(L)=0$. Hereafter $Z(A)$ stands for the center of a given algebra $A$.

So, it is natural to introduce  the notion of eigenvalues of the generating matrix $L$ as  roots of the characteristic polynomial.

\begin{remark} \rm Note that a similar CH identity is valid for the generating matrix of the RE algebra (modified or not), associated with any involutive or Hecke symmetry $R$, provided $R$ is  skew-invertible and even symmetry. The reader is referred to \cite{GPS1} for details.
 \end{remark}

 In the quantum double $(A,B)$, which we are going to construct, the role of the $A$-module $B$ will be played by the algebra $U(gl(N)_h)$. As for the algebra $A$, it is a unital associative algebra, generated by  elements $\partial_i^j$ which are entries of the $N\times N$ generating matrix $D=\|\pa_i^j\|_{1\leq i,j \leq N}$. The algebra $A$ is commutative. This property, expressed in a matrix form, reads
\be
D_1 D_2=D_2 D_1.
\label{dva}
\ee

Besides, we define  {\em permutation relations} between generators of the algebras $A$ and $B$ by the rule
\be
D_1 L_2=L_2 D_1+P_{12}+h D_1P_{12}.
\label{tri}
\ee
The system (\ref{tri}) enables us to define a map $\sigma: A\ot B\to B\ot A$ by sending any  element $\pa_i^j\, l_k^m$
from the left hand side of  (\ref{tri}) to the corresponding element from the right hand side. By  consecutive applying this rule and by assuming that
$$
\sigma(1_A\ot b)=b\ot 1_A\quad \forall \, b\in B\quad \mathrm{and} \quad \sigma(a\ot 1_B)=1_B\ot a\,\quad\forall\, a\in A,
$$
we can send any element
from $A\ot B$ to an element from $B\ot A$.  We call $\sigma$ {\em the   permutation map}. 

Note that the above algebras $A$ and $B$  are quotient-algebras:
$$
A=T(\span(\pa_i^j))/\langle J_1\rangle,\quad B=T(\span(l_i^j))/\langle J_2\rangle,
$$
where the notation $T(V)=\oplus_{k\ge 0} V^{\ot k}$ stands for the free tensor algebra, generated by a space $V$ and $\langle J_i\rangle$, $i=1,2$ are some ideals of $T(V)$. 
We say that the permutation relations are  compatible with the defining relations of the algebras $A$ and $B$ if the map $\sigma$  preserves these ideals, i.e.
$$
\sigma(J_1\ot \span(l_i^j))\subset \span(l_i^j)\ot J_1\quad \mathrm{and} \quad \sigma(\span(\pa_i^j)\ot J_2)\subset  J_2\ot \span(\pa_i^j).
$$

\begin{proposition}
The permutation map $\sigma$ is compatible with the systems of relations {\rm (\ref{odin})} and {\rm (\ref{dva})}.
\end{proposition}

A verification of the compatibility of the systems  (\ref{odin}), (\ref{dva}) and the permutation relations (\ref{tri})  is straightforward: it suffices to check that
 the following relation
$$
 D_1(L_2 L_3-L_3 L_2-hL_2P_{23}+h P_{23}L_2)=(L_2L_3-L_3 L_2-hL_2P_{23}+h P_{23} L_2) D_1
$$
is valid in virtue of the permutation relations.

\begin{remark} \rm
There exists another system of  permutation relations compatible with systems (\ref{odin}) and (\ref{dva}), namely, the following one
\be
D_1 L_2=L_2 D_1+P_{12}-hP_{12} D_1. \label{trii} \ee
However, the system (\ref{odin}), (\ref{dva}), and  (\ref{tri}) is equivalent to that (\ref{odin}), (\ref{dva}), and (\ref{trii}). In order to show  this equivalence
it suffices to apply the operation of transposition to all matrices entering the first system and to replace  $h$ by $-h$.
\end{remark}

Given a quantum double $(A,B)$, let us assume that there exists a counit (an algebra homomorphism) $\varepsilon_A:\,A\to \C$. Then we are able to construct an 
action  of the algebra $A$ onto $B$ by setting
\be
a\triangleright b=(I\ot \varepsilon_A)\sigma(a\ot b),\qquad a\in A,\, b\in B.
\label{gen-act}
\ee
Here, we identify $b\ot 1_{\C}$ and $b$.

Let us go back to the above quantum double. We define the counit  $\varepsilon_A$ in the algebra $A$ by setting
$$
\varepsilon_A(1_A)=1_{\C},\quad \varepsilon_A(\pa_i^j)=0\quad  \forall \, i,j,\quad \varepsilon_A(ab)= \varepsilon_A(a)\varepsilon_A(b).
$$

Then  formula (\ref{gen-act}) leads to the action $\pa_i^j\triangleright 1_B=0$, $\forall i,j$. Also, by taking into account (\ref{tri}), we find
$$
\pa_i^j\triangleright l_k^s=\de_i^s \de_k^j \quad   \Leftrightarrow \quad  D_1\triangleright L_2=P_{12}.
$$
 This action coincides with that of the usual partial derivatives in the generators of the commutative algebra $Sym(gl(N))$.

Now, we apply the QPD to a second order monomial. We have
\be
D_1\triangleright(L_2L_3) = P_{12}L_3+L_2P_{13} + hP_{12}P_{23}.
\label{44-1}
\ee
A general formula of applying the QPD to an arbitrary monomial can be found in \cite{GS1}.

Let us observe that there exists  a coproduct $\De:A\to A\ot A$, enabling us to express the above action as follows. We put $\De(1_A)=1_A\ot 1_A$ and
\be
\De(\pa_i^j)=\pa_i^j\ot 1_A+1_A\ot \pa_i^j +h\sum_k \pa_k^j\ot \pa^k_i:={\pa_i^j}_{\!(1)}\ot {\pa_i^j}_{\!(2)}.
\label{first}
\ee
Here, in the last equality we use the standard Sweedler's notation for the coproduct $\Delta(\partial_i^j)$. On the whole algebra $A$ the map $\Delta$ is extended by
the homomorphism property
$$
\De(\pa_i^j\, \pa_m^n)={\pa_i^j}_{\!(1)}\,{\pa_m^n}_{\!(1)}    \ot {\pa_i^j}_{\!(2)}\,{\pa_m^n}_{\!(2)}
$$
and so on.

\begin{proposition}
This  coproduct is coassociative. Also, the QPD applied according to the usual rule
$$
\pa_i^j\triangleright (a\, b)= ({\pa_i^j}_{\!(1)}\triangleright a)\, ({\pa_i^j}_{\!(2)}\triangleright b),
$$
send the elements \rm
\be
 l_k^s\, l_p^q-l_p^q\,l_k^s - h(\delta_p^sl_k^q - \delta_k^ql_p^s)
\label{55}
\ee
\it
to zero.
\end{proposition}

As follows from this proposition, the ideal $\langle J_2\rangle$ is invariant subspace with respect to the action of the QPD $\pa_i^j$.
Otherwise stated, the QPD extended via the above coproduct are well-defined on the algebra $U(gl(N)_h)$. Note that the coproduct $\De$ together with the counit $\varepsilon_A$
define a bi-algebra structure in the algebra $A$. 

Let us perform a linear shift of the generating matrix $D=\|\partial_i^j\|$ of the algebra $A$ as follows
\be
\hat{D}=D+\frac{I}{h} \quad \Leftrightarrow \quad \hpa_i^j=\pa_i^j+\frac{\delta_i^j}{h}\,1_A.
\label{modi}
\ee
The permutation relations (\ref{tri}) being expressed in terms of the matrix $\hD$ become
\be
\hD_1\, L_2=L_2\, \hD_1+h \hD_1P_{12}.
\label{triii}
\ee

The action of the copproduct $\Delta$ (\ref{first}) on the generators $\hpa_i^j$ reads
\be
\De(\hpa_i^j)=h\sum_k \hpa_k^j\ot \hpa_i^k \quad \Leftrightarrow \quad
\De(\hD^t)=h \hD^t\stackrel{.}{\ot} \hD^t,
\ee
where  $t$ stands for  the matrix transposition, and notation $M\stackrel{.}{\ot} N$ denotes the matrix with entries $m_i^k\ot n_k^j$ for any two square matrices
$M=\|m_i^j\|$ and  $N=\|n_i^j\|$ of the same size.

This coproduct leads to the following form of the Leibniz rule
\be
\hD(ab)^t=h \hD(a)^t\hD(b)^t,
\label{Lei}
\ee
where $\hD(a)$ denotes  the matrix with entries $\hpa_i^j(a)$. Below, we use the matrix $\DD=h \thD$, which is multiplicative (group-like) with respect to the coproduct:
\be
\De(\DD)=\DD\stackrel{.}{\ot} \DD.
\label{Leiii}
\ee

In conclusion, we want to emphasize that the method of defining analogs of partial derivatives  by permutation relations is valid  for the  modified RE algebra,
associated with any skew-invertible symmetry (involutive or Hecke), whereas the corresponding coproduct can be constructed only if $R$ is involutive (i.e. $R^2=I$).  Thus, the first method is more general.

\section{Quantum partial derivatives on $U(u(2)_h)$}

Let us consider the case $N=2$ in more detail. Below, we use the following notations $a=l_1^1,\, b=l_1^2,\, c=l_2^1,\, d=l_d^2$.
Thus, we have
$$
L = = \left(\!\!
\begin{array}{cc}
l_1^1&l_i^2\\
l_2^1&l_1^2
\end{array}
\!\!\right)=\left(\!\!
\begin{array}{cc}
a&b\\
c&d
\end{array}
\!\!\right),\qquad
\hD  = \left(\!\!
\begin{array}{cc}
\hat \partial_1^1&\partial^2_1\\
\partial_2^1&\hat\partial_2^2
\end{array}
\!\!\right)= \left(\!\!
\begin{array}{cc}
\hat \partial_a&\partial_c\\
\partial_b&\hat\partial_d
\end{array}
\!\!\right).
$$

Hereafter, we use  the ``hat'' notation only for the diagonal elements of the matrix $\hD$ since the off-diagonal partial derivatives are actually not changed under the shift (\ref{modi}):
$\hat \partial_i^j\equiv \partial_i^j$ for any $i\not= j$.

Let us exhibit the permutation relations (\ref{triii}) in the case under consideration
$$
\begin{array}{llllll}
\hpa_a \,a = a\,\hpa_a+h \, \hpa_a & \hspace*{2mm} &\hpa_a\,b = b\, \hpa_a + h\,\pa_c& \qquad \pa_c\,a = a\, \pa_c& \hspace*{2mm} & \pa_c\,b = b\, \pa_c\\
\rule{0pt}{4.5mm}
\hpa_a\, c = c\,\hpa_a & \hspace*{2mm} & \hpa_a\,d = d\, \hpa_a &\qquad \pa_c\,c = c\, \pa_c + h\, \hpa_a & \hspace*{2mm} & \pa_c\,d = d\, \pa_c + h\,\pa_c\\
\rule{0pt}{6mm}
\pa_b \,a = a\, \pa_b + h\,\pa_b & \hspace*{2mm} & \pa_b\,b = b\, \pa_b + h\,\hpa_d&\qquad \hpa_d\,a = a\, \hpa_d & \hspace*{2mm} & \hpa_d\,b = b\, \hpa_d\\
\rule{0pt}{4.5mm}
\pa_b \,c = c\,\pa_b & \hspace*{2mm} & \pa_b\,d = d\, \pa_b & \qquad \hpa_d\,c = c\, \hpa_d + h\,\pa_b & \hspace*{2mm} & \hpa_d\,d = d\, \hpa_d + h\,\hpa_d.
\end{array}
$$

Now, we pass  to the compact form $u(2)_\h $ of the algebra $gl(2)_\h$. Namely, we introduce the following generators
$$
t={{1}\over{2}}(a+d),\qquad     x={{i}\over{2}}(b+c),\qquad y= {{1}\over{2}}(c-b),\qquad
z={{i}
\over{2}}(a-d).
$$
The corresponding Lie brackets read:
$$
[x, \, y]=\h z,\qquad [y, \, z]=\h x,\qquad[z, \, x]=\h y,\qquad [t, \, x]=[t, \, y]=[t, \, z]=0.
$$

Since the change of the generators is linear, we apply the classical formula for computing the QPD in new generators: $\partial_t = \pa_a+\pa_d$ and so on. In these generators the
"shifted version" is used only for the derivative in $t$:  $\hpa_t = \pa_t + \frac{2}{\h}\,{\rm id}$.
The permutation relations  become
\be
\begin{array}{l@{\quad}l@{\quad}l@{\quad}l}
\hpa_t\,t - t\,\hpa_t = \frac{h}{2}\,\hpa_t & \hpa_t\, x - x\,\hpa_t
=-\frac{h}{2}\,\pa_x &
\hpa_t\, y - y\, \hpa_t=-\frac{h}{2}\,\pa_y &\hpa_t\, z - z\,\hpa_t=- \frac{h}{2}\,\pa_z\\
\rule{0pt}{7mm}
\pa_x\, t - t\,\pa_x = \frac{h}{2}\,\pa_x &\pa_x \,x -  x\,\pa_x = \frac{h}{2}\,\hpa_t &
\pa_x \, y-  y\,\pa_x = \frac{h}{2}\,\pa_z & \pa_x \,z - z\, \pa_x  = - \frac{h}{2}\,\pa_y \\
\rule{0pt}{7mm}
\pa_y \,t - t \, \pa_y = \frac{h}{2}\,\pa_y & \pa_y \,x -  x\,  \pa_y = -\frac{h}{2}\,\pa_z &
\pa_y \,y - y \,  \pa_y = \frac{h}{2}\,\hpa_t & \pa_y \,z - z \,  \pa_ y= \frac{h}{2}\,\pa_x\\
\rule{0pt}{7mm}
\pa_z \,t - t \,\pa_z = \frac{h}{2}\,\pa_z & \pa_z \,x - x \,\pa_z = \frac{h}{2}\,\pa_y&
\pa_z \,y -  y\,\pa_z = -\frac{h}{2}\,\pa_x & \pa_z \,z - z \,\pa_z = \frac{h}{2}\,\hpa_t.
\end{array}
\label{leib}
\ee

These  permutation relations play the role of the Leibniz rule for the QPD acting on the algebra $U(u(2))$. Another form of this rule can be expressed via a coproduct acting on the QPD. 
In a matrix form this coproduct reads
\be 
\De(\hatt)=\hatt\stackrel{.}{\ot}\hatt,\qquad
{\hatt}=i\hh \left(\begin{array}{rrrr}
\hpa_t&\pa_x&\pa_y&\pa_z\\
-\pa_x&\hpa_t&-\pa_z&\pa_y\\
-\pa_y&\pa_z&\hpa_t&-\pa_x\\
-\pa_z&-\pa_y&\pa_x&\hpa_t
\end{array} \right),
\label{seven}
\ee
Hereafter, we use  a new parameter $\hh=h/2i$. Also,  the symbol $\hatt(a)$  denotes the matrix whose entries  result from applying  the corresponding partial derivatives to an element
$a\in U(u(2)_h)$. Thus, we get a linear map $U(u(2)_h)\to \mathrm{Mat}_2(U(u(2)_h))$ which  will be also denoted by the  symbol $\hatt$:
\be
\hatt: \, a\mapsto  \hatt(a)\quad \, \forall\, a\in U(u(2)_h).
\label{mapp}
\ee

We state that the map $\hatt$  preserves the multiplication of elements of the algebra $U(gl(2)_h)$
\be
\hatt(ab)=\hatt(a) \hatt(b),\quad \forall\,a,b\in U(u(2)_h)
\label{mult}
\ee
 Thus, the matrix $\hatt$ enables us to present the Leibniz rule in a multiplicative or group-like form.
 Consequently, the map $\hatt$ defines a  representation of the algebra $U(u(2)_h)$:
 $$
\hatt(x)\hatt(y)-\hatt(y)\hatt(x)=2i\hh\hatt(z)=h\hatt(z)
$$
and so one.

The images of the $U(u(2)_h)$ generators under the map $\hatt$ are as follows:
\be
\begin{array}{lcl}
\hatt(1_{U(u(2)_h)})=I,&\quad &\hatt(t)=(t+i\hh) I,\\
\rule{0pt}{5mm}
\hatt(x)=x\, I+i\hh  A,&\quad &\hatt(y)=y\, I+i\hh  B,\qquad \hatt(z)=z\, I+i\hh\, C,
\end{array}
 \label{mat0}
\ee
where
\be
A=
\left(\!\!
\begin{array}{cccc}
 0 &    1&0&0\\
 -1&0&0&0\\
0&0&  0 & -1\\
0&0& 1 & 0
\end{array}
\!\!\right),\qquad
B=
\left(\!\!
\begin{array}{cccc}
 0 & 0\,&1&0\\
0 &0&0&1\\
-1&0&  0 & 0\\
0&-1&  0 & 0
\end{array}
\!\!\right),\qquad
C=
\left(\!\!
\begin{array}{cccc}
 0 & 0&0&1\\
 0 &0&-1&0\\
0&1&  0& 0\\
-1&0& 0 & 0
\end{array}
\!\!\right).
\label{ABC}
\ee

The numerical matrices $A$, $B$ and $C$ possess the following multiplication table
\be
A^2=B^2=C^2=-I,\quad AB=-BA=C,\quad BC=-CB=A,\quad CA=-AC=B.
\label{mt}
\ee

Below, we extend the action of the QPD onto some lager  algebras with preserving the property  (\ref{mult}) of the map $\hatt$.

\section{QPD on central extension of $U(u(2)_h)$}

First, we present the CH identity for the generating matrix of the algebra $U(gl(2)_h)$ expressed in terms of the compact generators
$$
L=\left(\!\!
\begin{array}{cc}
 a & b\\
 c &d
\end{array}\!\!
\right)=\left(\!\!
\begin{array}{cc}
 t-iz & -ix-y\\
 -ix+y& t+iz
\end{array}\!\!
\right)
$$
This identity reads $p(L)=0$, where the characteristic polynomial $p(\tau)$ is of the form
$$
p(\tau)=\tau^2-(2t+h)\,\tau+(t^2+x^2+y^2+z^2+2i\hh\,t)\, I.
$$
Let us denote $\mu_1$ and $\mu_2$ the roots of the characteristic polynomial $p(\tau)$:
$$
\mu_1+\mu_2=2t+2i \hh,\qquad \mu_1\, \mu_2=t^2+\Cas+2i\hh \,t,
$$
where we use the notation $\Cas=x^2+y^2+z^2$. So, the elements $\mu_i$ belong to a central extension of  $U(u(2))_h$.

Consider the algebra $\AA=U(u(2)_h)[\mu_1, \mu_2]$ and compute the quantities $\hatt(\mu_i),\,\,i=1,2$. With the use of (\ref{mat0}) it is not difficult to
find their sums:
\begin{eqnarray}
&&\hatt(\mu_1)+\hatt(\mu_2)=\hatt(\mu_1+\mu_2)=\hatt(2t+2i\hh) = (2t+4i\hh)I,
\label{razz}
\end{eqnarray}

Now, we  find the difference $\hatt(\mu_1)-\hatt(\mu_2)$. Consequently, it becomes possible to calculate each $\hatt(\mu_i)$. For this purpose let us first compute the matrix $\hatt(\mu^2)$,
where $\mu=\mu_1-\mu_2$. Taking into account that
$$
\mu^2=(\mu_1+\mu_2)^2-4\mu_1\, \mu_2=-4(\Cas+\hh^2),
$$
and  with the use of (\ref{mat0}) we have
$$
\hatt(\mu^2) = \hatt(-4(x^2+y^2+z^2+4\hh^2)) =(\mu^2 -12\hh^2)I-8i\hh\,M,
$$
where $M= (xA+yB+z\,C)$.

Now, to find the matrix $\hatt(\mu)$ we have to calculate the square root $\hatt(\mu)=\sqrt{\hatt(\mu^2)}$.
This relation is motivated by our wish to preserve the multiplicative property of the map $\hatt$.

Since the matrix $M=xA+yB+z\,C $ satisfies the CH identity
$$
M^2-2i\,\hh\, M+\Cas\, I=0
$$
its eigenvalues are easy to find:
$$
\la_1=i \hh+\frac{\mu}{2},\qquad \la_2=i\hh-\frac{\mu}{2}.
$$
Consequently, the eigenvalues of the matrix $\hatt(\mu^2)$ are as follows
$$
\nu_1=\mu^2-4\hh^2-4i\hh\,\mu=(\mu-2i\hh)^2,\qquad \nu_2=\mu^2-4\hh^2+4i\hh\,\mu=(\mu+2i\hh)^2.
$$

As for the matrix $\hatt(\mu)$, it can be found via the spectral decomposition of the matrix $\hatt(\mu^2)$:
$$
\hatt(\mu)=\frac{\hatt(\mu^2)-\nu_2\, I}{\nu_1-\nu_2}\,\sqrt{\nu_1} +\frac{\hatt(\mu^2)-\nu_1\, I}{\nu_2-\nu_1}\,\sqrt{\nu_2}.
$$
Thus, we have 4 candidates to the role of the matrix $\hatt(\mu)$ in dependence of the sign choices for the roots
$$
\sqrt{\nu_1}=\epsilon_1(\mu-2i\hh),\qquad \sqrt{\nu_2}=\epsilon_2(\mu+2i\hh),\qquad \epsilon_1, \epsilon_2 \in \{\pm 1\}.
$$

But a straightforward verification shows that the only choice $\epsilon_1=\epsilon_2=1$ leads to the result compatible with the classical limit (i.e. corresponding to $\hh=0$) for the
action of the partial derivatives:
\be
\pa_t \mu=0\quad {\rm  and}\quad \pa_x\, \mu=-4\,\frac{x}{\mu}.
\label{cla}
\ee
The second relation in (\ref{cla}) follows from the equality $\mu^2=-4\,\Cas$ valid in the classical limit.

Thus, we arrive at the final formula
\be
\hatt(\mu)=\frac{(\mu+2i\hh)^2+(\mu-2i\hh)^2}{2\mu}\, I+\frac{(\mu-2i\hh)-(\mu+2i\hh)}{\mu}\,M=\frac{\mu^2-4\hh^2}{\mu}\, I-\frac{4i\hh}{\mu}\, M.
\label{sto}
\ee
In a similar way we can get the action of the QPD on an arbitrary power of $\mu$:
$$
\hatt(\mu^p)=\frac{(\mu+2i\hh)^{p+1}+(\mu-2i\hh)^{p+1}}{2\mu}\, I+\frac{(\mu-2i\hh)^{p}-(\mu+2i\hh)^{p} }{\mu}\, M.
$$

Now, by using  (\ref{razz}) and (\ref{sto}) we get
$$
\hatt(\mu_1)=\left(t+\frac{(\mu^2+2i\hh)^2}{2\mu}\right) I-\frac{2i\hh}{\mu}\, M,\qquad \hatt(\mu_2)=\left(t-\frac{(\mu - 2i\hh)^2}{2\mu}\right) I+\frac{2i\hh}{\mu}\, M.
$$
These formulae entail, for example, the following
$$
\hpa_t(\mu_1)=-\,\frac{i}{\hh}\,\left(t+\frac{(\mu+2i\hh)^2}{2\mu}\right ),\quad \pa_x(\mu_1)=-2\,\frac{x}{\mu},
$$
and so on.

Thus, we have computed the action of the QPD onto the elements $\mu_i$, $i=1,2$. We complete this section by recalling the notion of the quantum radius and by finding
the result of QPD action on this element.

In our previous papers we introduced the so-called quantum radius $\rh=\sqrt{\Cas+\hh^2}$. In terms of $\mu$ it is written as $\rh = \pm\mu/ 2i$. The sign here is not fixed because our ordering of
$\mu_1$ and $\mu_2$ is arbitrary.

By using the above relation between $\mu$ and $\rh$ (with any sign) we get from (\ref{sto})
\be
{\hatt}(\rh)=\frac{\rh^2+\hh^2}{\rh}\, I+\frac{i\hh}{\rh}\, M.
\label{mat1}
\ee
This formula enables us to find the action of all QPD on the quantum radius:
\be
\pa_t \rh=-\frac{i\hh}{\rh},\qquad \pa_x\rh=\frac{x}{\rh},\qquad \pa_y\rh=\frac{y}{\rh},\qquad \pa_z\rh=\frac{z}{\rh}.
\label{pervy}
\ee
Note that in the classical limit $\hh=0$ the derivative $\pa_t r$ vanishes and other formulae turn into similar ones (with $\rh$ replaced by $r$).

Now, consider the algebra
\be
\AA=\left(U(su(2)_h)\ot \C[t,\rh]\right)/\langle x^2+y^2+z^2+\hh^2-\rh^2 \rangle.
\label{quot}
\ee
It is easy to see that the map $\hatt$ sends the element $x^2+y^2+z^2+\hh^2-\rh^2$ to 0.  Consequently, this map is well-defined on the algebra $\AA$.

We assume  the quantum radius to be real and positive provided that the parameter $\hh$ is real and the generators $x$, $y$ and $z$ are represented by Hermitian operators.
Observe that if the algebra  $U(u(2)_h)$ is represented in the space of spin $n$ the quantum radius takes the value $\rh=(2n+1)\hh$.

Recall that in our treatment of the algebra $U(u(2)_h)$ as a non-commutative analog of the polynomial algebra on the Minkowski space, the generators $x,y, z$ play the role of  spacial variables and $t$ is interpreted as the time (see \cite{GS2}).

\section{Extension of  QPD on some elements of $\BB$ via CH identities}

In this section we extend the QPD onto some elements of  the skew-field ${\BB}=\AA[\AA^{-1}]$ with preserving
the Leibniz rule in its multiplicative form. Note that this Leibniz rule can be expressed via the matrix  $\hatt$ or $\DD$, introduced at the end of the section 2. 
Since the size of the matrix $\DD$ is smaller (this fact becomes more significant in the higher dimensions),
we deal with the matrix $\DD$ in this section. 

Let $b$ be an arbitrary nontrivial element of the algebra $\AA$. If we can extend the map $\DD$ onto the element  $b^{-1}\in \BB$ with preserving the Leibniz rule, we
should have
$$
\DD(b^{-1})=\DD(b)^{-1}.
$$
Thus, in order to compute  $\DD(b^{-1})$  we have to invert the matrix $\DD(b)$ with non-commutative entries. In principle, this procedure can be performed by means of 
the Gelfand-Retakh method using the so-called quasideterminants. In order to present the entries of the matrix $\DD(b^{-1})$ as elements of $\BB$ we need the Ore property of the algebra $\AA$. Hopefully, the algebra $\AA$ has this property since it is so for the algebra $U(u(2)_h)$. Nevertheless, practically, reduction of any "left"
fraction to a "right" one and vice versa in this algebra is difficult indeed. 

Below, we deal with some elements $b\in \AA$ for which the computation of the matrix $\DD(b)^{-1}$ can be performed by means of the CH identities for matrices with non-commutative  entries. 
For these elements $b$ we succeeded in finding the matrices $\DD(b)^{-1}$ with entries from the skew-field $\BB$. 

First, we calculate the matrix $\DD(b)$ for some basic $b\in \AA$. Taking into account the explicit form of $\DD$
\be
\DD=i\, \hh\, \left(\!\!
\begin{array}{cc}
 \hpa_t+i\pa_z &i\pa_x-\pa_y\\
\rule{0pt}{5mm}
i\pa_x+\pa_y&\hpa_t-i\pa_z
\end{array}
\!\!\right),
\label{D-expl}
\ee
we find
$$
\DD(t)=(t + i\hh)I,\quad
\DD(x)=\left(\!\!
\begin{array}{cc}
x&-\hh\\
-\hh&x
\end{array}
\!\!\right),\quad
\DD(y)=\left(\!\!
\begin{array}{cc}
y&-i\hh\\
i\hh&y
\end{array}
\!\!\right),\quad
\DD(z)=\left(\!\!
\begin{array}{cc}
z-\hh&0\\
0&z+\hh
\end{array}
\!\!\right)
$$
$$
\DD(\rh)=\frac{1}{\rh}\left(\!\!
\begin{array}{cc}
\rh^2+\hh^2-\hh\,z&-\hh\,(x+iy)\\
\rule{0pt}{5mm}
-\hh\,(x-iy)&\rh^2+\hh^2+\hh\,z
\end{array}
\!\!\right).
$$
The  matrix $\DD(\rh)$ obeys the CH identity:
\be
\DD(\rh)^2 - 2\rh\,\DD(\rh) + (\rh^2 - \hh^2) I = 0,
\label{char}
\ee
which can be easily rewritten in the factorized form:
\be
(\DD(\rh) - (\rh+\hh)I)(\DD(\rh)-(\rh-\hh)I) =0.
\label{d-ch-f)}
\ee

Now, consider the following example. Let $b$ be a linear combination of generators with  coefficients $\alpha_i \in \C$:
\be
b=\al_0\, t+\al_1\, x+ \al_2\, y+\al_3\, z,\qquad  \al_1^2+\al_2^2+\al_3^2 \not= 0.
\label{bbb}
\ee
Since
$$
\DD(b)=(b+i\al_0)\, I-\hh\, N,\qquad
N=\left(\!\!
\begin{array}{cc}
\al_3&\al_1+i\al_2\\
\al_1-i\al_2&-\al_3
\end{array}
\!\!\right),
$$
then taking into account  that $N^2=(\al_1^2+\al_2^2+\al_3^2)I$, we find the matrix  $\DD^{-1}(b)$:
$$
\DD(b)^{-1}=\frac{(b+i\al_0)I+\hh N}{(b+i\al_0)^2-(\al_1^2+\al_2^2+\al_3^2)\hh^2}.
$$
Note that the meaning of the fraction here is not ambiguous since its numerator and denominator commute with each other. 

Introduce an element $c=\rh-b$, where $b$ is given by (\ref{bbb}) with the following restriction on the coefficients
$$
\al_0=0,\qquad \al_1^2+\al_2^2+\al_3^2 =1.
$$
The matrices $\DD(b)$ and $\DD(\rh)$ commute with each other and obey the second order CH identities 
\be
(\DD(b)-b\, I)^2= (\al_1^2+\al_2^2+\al_3^2)\, \hh^2\, I,
\label{char-db}
\ee
and  (\ref{char}) respectively. Thus, it is reasonable to look for the inverse
matrix $\DD(c)^{-1}$ of the form:
\be
 \DD(c)^{-1}=I\, a_0+\DD(\rh)\, a_1+\DD(b)\,a_2+ \DD(\rh)\DD(b)\, a_3.
\label{last}
\ee

On multiplying (\ref{last}) by the matrix $\DD(c)=\DD(\rh)-\DD(b)$ from the left and demanding the result to be the unit matrix $I$ we find a system of linear equations
for coefficients $a_i$:
\begin{eqnarray*}
-(\rh^2-\hh^2)\, a_1+(b^2-\hh^2)\, a_2=1&&\\
a_0+2\rh\, a_1+(b^2-\hh^2)\, a_3=0&&\\
a_0+2b\, a_2+(\rh^2-\hh^2)\,a_3=0&&\\
a_1-a_2-2(\rh-b)\, a_3=0
\end{eqnarray*}
Here, we used the relations (\ref{char-db}) and (\ref{char}). 

The coefficients of the system above commute with each other, so the solution can be found by the standard methods of linear algebra. 
Applying the Cramer's rules we get the following result
$$
a_0=\frac{2(\rh^2+b^2-3b\rh-\hh^2)}{(\rh-b)((\rh-b)^2-4\hh^2)},\qquad a_1=\frac{3b-\rh}{(\rh-b)((\rh-b)^2-4\hh^2)},
$$
$$
a_2=\frac{3\rh-b}{(r-b)((\rh-b)^2-4\hh^2)},\qquad a_3=\frac{-2}{(\rh-b)((\rh-b)^2-4\hh^2)}.
$$

Now, we are able to extend the action of the QPD on the element $c^{-1} = (\rh-b)^{-1}$. Namely, from  (\ref{D-expl})  we get
$$
(\DD_1^1+\DD_2^2)(c^{-1}) = 2i\hh\,\hat \partial_t c^{-1},\qquad (\DD_1^2+\DD_2^1)(c^{-1}) = -2\hh\, \partial_x c^{-1},\quad \mathrm{etc.}
$$
To write down the general answer it is convenient to introduce the following vectors
$$
\vec \alpha = (\alpha_1,\alpha_2,\alpha_3),\qquad \vec \rh = (x,y,z),\qquad \vec \nabla =(\partial_x,\partial_y,\partial_z).
$$

Then, in virtue of the relation $\DD(c^{-1}) = \DD(c)^{-1}$ we come to the final results:
$$
\partial_t\left(\frac{1}{\rh-b}\right) = \frac{-i\hh}{\rh((\rh-b)^2-4\hh^2)},
$$
$$
\vec\nabla\left(\frac{1}{\rh-b}\right) = \frac{\rh\vec\alpha - \vec \rh}{\rh}\,\frac{1}{((\rh-b)^2-4\hh^2)} - (\hh\,\vec\alpha +i[\vec \rh\times\vec\alpha])\,\frac{2\hh}{\rh(\rh-b)((\rh-b)^2-4\hh^2)},
$$
where $[\cdot \times\cdot]$ stands for the vector product of two vectors. We point out that the components of the vector $\vec \rh$ do {\it not} commute with the
element $\rh-b$ entering the denominator. So, we have to explicitly fix the order of factors in the formula for $\vec \nabla (\rh-b)^{-1}$.

Also, note that in the limit $\hh\to 0$ we recover the classical result:
$$
\vec \nabla\left(\frac{1}{r-b}\right)=\frac{r\vec\alpha - \vec r}{r(r-b)^2},
$$
where  $r=|\vec r|$ and  $b= \vec\alpha\cdot\vec r$ is the scalar product. 

Thus, by applying  the QPD to  elements of the form $a c^k$, $a\in \AA$, $k\in {\Bbb{Z}}$, we can represent the resulting  elements as those from $\BB$.

Our interest in applying the QPD to the element $(\rh-b)^{-1}$ is motivated by our wish to construct
a non-commutative analog of the Dirac potential. Let us recall that in the  classical setting this potential is a solution of the following
equation ${\rm rot}\, \AAA=\HH$, where $\HH$ is the vector of the magnetic field, namely, a stationary solution of the Maxwell system 
\be
 {\rm div}\, \HH=0,\qquad {\rm rot}\, \HH=4 g \pi \de(r),
\label{Dm} 
\ee
where $g$ is a nontrivial constant factor.
The electric field is assumed to vanish. 
As usual, the notation ${\rm div}$ and ${\rm rot}$ stand for respectively the divergence and rotor of a given vector field, and $\de(r)$ is the delta-function.
Dirac found a solution to the  system (\ref{Dm}) in the form $\HH=g\,\frac{{\rm \bf r}}{r^3}$.

Besides, he found a family of vector-potentials  of this model:
\be
\AAA=\frac{g}{r}\frac{[\vec{r}\times \vec \al]}{(r-\vec{r}\cdot \vec{\al})}, 
\label{pot}
\ee
where, $\vec{\al}$ is a unit vector,  and $\cdot$ is the scalar product of two vectors. Each of these vector-potentials
is singular on a half-line.

In \cite{GS2} we have found a non-commutative counterpart of the Dirac monopole, i.e. a solution to the system (\ref{Dm}). It is of the form
\be
\H=\frac{g\, \vec{\rh}}{\rh (\rh^2-\hh^2)} .
\label{mag}
\ee

Unfortunately, we have not succeeded in finding a non-commutative counterpart of the potential $\AAA$, i.e. a solution to the equation ${\rm rot}\, \AAA=\HH$, where $\HH$ is defined
by (\ref{mag}). It should be emphasized that the problem of finding such a potential  in the framework of our non-commutative setting is much more complicated, than in the classical setting.
We want only to note that by  looking for a solution of the equation  ${\rm rot}\, \AAA=\HH$    under a form similar to the classical one, we have to take in consideration  the fact that the element $(\rh-b)^{-1}$ does not commute with the components of the vector $\vec{\rh}$.

\end{document}